\documentclass[12pt]{article}
\usepackage{amssymb}

\title{Heat Kernels, Symplectic Geometry, Moduli Spaces and Finite Groups}

\author{Kefeng Liu}

\date{}

\begin{document}

\maketitle


\section{Introduction} 

In this note we want to discuss some applications of heat kernels in symplectic geometry, moduli spaces and finite groups. More precisely we will prove the nonabelian localization formula in symplectic geometry, derive formulas for the symplectic volume and intersection numbers of the moduli space of flat connections on a Riemann surface, and obtain several quite general formulas for the numbers of solutions of equations in finite groups.  Several new formulas for the push-forward measures by various maps between Lie groups are also obtained. In solving these problems, we will use the corresponding heat kernels on Euclidean spaces, on Lie groups and on finite groups. 

The discussions of some aspects of the first two applications have appeared in [Liu], [Liu1] and [Liu2], here we will only sketch the main ideas. The purpose to include them here is to unify the discussions by using heat kernels. The third and the fuorth applications were worked out through many discussions with P. Diaconis. Several results for compact Lie groups in Section 4 were also motivated by his conjectures and his results for finite groups 

The main idea we use is very simple and goes as follows. We consider a map between two spaces 
$f: \ M\rightarrow N$. In many problems we are interested in understanding the inverse image $f^{-1}(x_0)$ for some point $x_0\in N$. Assume that there is a heat kernel $H(t, x, x_0)$ on $N$. We consider the integral of the pull-back of $H(t, x, x_0)$ by $f$:

$$I(t)= \int_M H(t, f(y), x_0)d\,y$$
with respect to certain measure $dy$ on $M$. We then perform the computations of this integral in two different ways, similar to the heat kernel proofs of the Atiyah-Singer index formula and the Atiyah-Bott fixed point formula. On one hand,  when $t$ goes to $0$, $I(t)$ will localize to an integral over a neighborhood of $f^{-1}(x_0) \subset M$, on the other hand, by using the special properties of symplectic manifolds or representation theory of the corresponding groups, we can explicitly compute $I(t)$ globally in our cases. In this way, we obtain formulas about the geometry and topology of $f^{-1}(x_0)$ in terms of certain global information on $M$ and $N$. We will give several  examples here,  to all of which  we can apply this simple idea.

(1). The moment map $\mu:\ \ M\rightarrow g^*$, where $M$ is a symplectic manifold, and $g^*$ is the dual of the Lie algebra of the group $G$ acting on $M$. In this case we are interested in the symplectic reduction $M_G=\mu^{-1}(0)/G$.

(2). The holonomy map $f:\ \ G^{2g}\rightarrow G$ where $G$ is a compact Lie group and $f$ is the product of the commutators 

$$f(x_1, y_1, \cdots x_g, y_g)=\prod^g_{j=1}[x_j, y_j].$$ In this case we are interested in the space $f^{-1}(c)/{ Z}_c$, which is the moduli space of flat connections on a flat $G$-bundle on a Riemann surface. Here $c\in G$ and $Z_c$ denotes the stabilizer of $c$. More general examples were studied in [Liu2].  

(3). Equivalently we can consider the map

 $$f:\ \ G^{2g}\times \prod_{j=1}^n {O}_{c_j} \rightarrow G$$ 
where ${O}_{c_j}$ is the conjugacy class in $G$ passing through $c_j \in G$. 

Another interesting map is the $n$-commutator map $f^{(n)}: \ \ G^n \rightarrow G$ with 

$$f^{(n)}(x_1, x_2, \cdots, x_n)= [x_1, [x_2, \cdots [x_{n-1}, x_n]]].$$ In this case, an inductive formula for the push-forward by $f^{(n)}$ of the Riemannian measure on $G^n$ can be obtained. 

(4). For a finite group $G$, we can consider the maps similar to those in (2) and (3) or more generally we can consider 

$$f: \ \ G^n \times \prod_{j=1}^n H_j \rightarrow G$$ where $H_1, \cdots, H_n$ are  subgroups of $G$ and 

$$f(x_1, \cdots, x_n; z_1, \cdots, z_n)=\prod_{j=1}^nx_j z_j x_j^{-1}.$$In this case the integral to define $I(t)$ is replaced by the sum over $G^n\times \prod_j H_j$, and the limit as $t$ goes to $0$ gives the number of solutions to the equation in $ G^n \times \prod_{j=1}^nH_j $: 

 	$$\prod_{j=1}^nx_j z_j x_j^{-1}=e.$$

More general examples than the above cases will also be considered in this note.

The interesting point here is that, the heat kernel method supplies a unified way  to deal with several seemingly unrelated problems in geometry, topology and finite group theory. In Section 2 we derive a nonabelian localization formula in symplectic geometry from the heat kernel point of view. Here we use the heat kernel of Euclidean space. Note that, different treatment of the result has been discussed first by Witten [W1], then in [Wu], [JK] and [Liu]. In Section 3 we obtain the formulas for the symplectic volume and some intersection numbers of the moduli space by using the heat kernel of the Lie group $G$. Some detailed discussion in this section has appeared in [Liu1], [Liu2], and in [BL]. In Section 4 we derive various formulas for the push-forward by those maps in (2) and (3) of the Riemaniann measures. In Section 5 we derive several formulas for counting the numbers of solutions of certain polynomial equations in finite groups. All of the results in this section grew out of discussions with P. Diaconis who has quite different proofs for these formulas by using combinatorics methods. One of the formulas in Section 5 has been proved and used as the main tools in [St] to solve the long-standing Brauer $p$-block conjecture. 

I want to dedicate this paper to the memory of Prof. Qi-Ming Wang who passed away ten years ago. In 1988, after reading my master thesis, Prof. Wang encouraged me to follow Prof. S.-T. Yau to study geometry. This is the turning point of my career.
   
Finally I would like to thank P. Diaconis for the many stimulating and enlightening discussions. 

\section{Symplectic reduction}

Let $M$ be a compact symplectic manifold with a symplectic form $\omega$. Assume that the compact Lie group $G$ acts on $M$. Let $g$ denote the Lie algebra of $G$ and $g^*$ be its dual. Assume that the $G$-action has a moment map 

$$\mu: \ \ M\rightarrow g^*.$$

With the metric induced from the Killing form on $G$, we can identify $g^*$ to the Euclidean space ${\bf  R}^n$. Let us denote the  metric on $g^*$ by $<\cdot, \cdot>$. Let 

$$H(t, x, x_0)=\frac{1}{(4\pi t)^{\frac{n}{2}}}{\mbox{ exp}}(-\frac{\|x-x_0\|^2
}{4t})$$ be the heat kernel on $g^*$. Assume $0\in g^*$ is a regular value of $\mu$. We are interested in studying the symplectic reduction $M_G=\mu^{-1}(0)/G$. We consider the integral of the pull-back of the heat kernel against the symplectic volume on $M$ :

$$ I(t)= \int_M H(t, \mu(y), 0)e^{\omega}.$$
As in the derivation of the index formula, we will compute this integral in two different ways: local and global. 

(i) Local computation. Let $t$ go to $0$, then the integral $I(t)$ localizes to a neighborhood of $\mu^{-1}(0)$, which can be identified as $\mu^{-1}(0)\times B_\delta$ where $B_\delta$ denotes a small ball of radius $\delta$ in $g^*$. More precisely we can write the above integral as 

$$I(t)= \int_{\mu^{-1}(0)\times B_\delta}H(t, \mu(y), 0)e^{\omega}+O(e^{-\delta^2/4t}).$$ 

When restricted to $\mu^{-1}(0)\times B_\delta$, the symplectic structure  $\omega$ has a canonical expression in terms of the induced symplectic form $\omega_0$ on the symplectic reduction $M_G=\mu^{-1}(0)/G$: 

$$\omega =\pi^* \omega_0+ d(\alpha, \theta)$$ where $\pi: \ \mu^{-1}(0)\rightarrow M_G$ denotes the quotient map which gives a principal $G$-bundle, and  $\theta$ is a connection form of this bundle. Also $\alpha$ denotes the coordinate on $g^*$ and $(\alpha, \theta)$ denotes the obvious paring. And in this local model $\mu$ is just the simple projection map:  $\mu: \ \mu^{-1}(0)\times g^*\rightarrow g^*$ defined by $\mu(x, \alpha)=\alpha$. 

By substituting the local expression of $\omega$ into the integrand we obtain  

$$I(t)=\int_{\mu^{-1}(0)\times g^*}H(t, \alpha, 0)e^{\pi^* \omega_0+ d(\alpha, \theta)}+O(e^{-\delta^2/4t}).$$

We now use the equalities: 

$$d(\alpha, \theta)=(d\alpha, \theta)+(\alpha, d\theta),\ {\mbox{and}}\  F=d\theta-\theta\wedge \theta$$
where $F$ denotes the curvature of $\theta$. From these we derive 

$$I(t)= \int_{\mu^{-1}(0)}e^{\pi^*\omega_0}\wedge \theta\int_{g^*}H(t, \alpha, 0)e^{(\alpha,F)} D\alpha +O(e^{-\delta^2/4t})$$ where $D\alpha$ denotes the volume form of $g^*$ and $\wedge \theta= \theta_1\wedge \cdots \wedge \theta_n$. Here we have used the fact that, by a degree count, in $e^{(d\alpha, \theta)}$ only the term $D\alpha \wedge \theta$ will contribute to the integral. Note that $\wedge \theta$ is actually a volume form of $G$.

By performing a simple Gaussian integral on $g^*\simeq {\bf{ R}}^n$, we get

$$I(t)= |G|\int_{M_G}e^{\omega_0-t<F, F>}+ O(e^{-\delta^2/4t})$$ where $<\cdot, \cdot>$ also denotes the inner product on $g$ induced from the Killing form, and $|G|$ denotes the corresponding volume of $G$.

(ii) Global computation. On the other hand we can rewrite $I(t)$ in terms of equivariant cohomology class.

$$I(t)= \int_g\int_M e^{-t<\varphi, \varphi>}e^{\omega+i(\mu, \varphi)}d\varphi$$

$$= \int_g e^{-t<\varphi, \varphi>}\int_M e^{\omega+i(\mu, \varphi)}d\varphi .$$Here $<\varphi, \varphi>$ denotes the inner product induced on $g$ by the Killing form. Note that the first identity actually corresponds to the Fourier transform of $H(t, x, 0)$. By comparing the computations in (i) and (ii), we have the following equality of Witten which we summarize as a proposition.

\vspace{.2in}

{\bf Proposition 1}: {\em We have the following formula: }

$$|G|\int_{M_G}e^{\omega_0-t<F, F>}={\mbox{lim}}_{t\rightarrow 0}\int_g\int_M e^{-t<\varphi, \varphi>}e^{\omega+i(\mu, \varphi)}d\varphi$$

\vspace{.1in}

This formula expresses the geometric and topological information on $M_G$ in terms of the information on $M$. This is the spirit of the nonabelian localization of Witten.

By using an idea of Wu [Wu], we can actually go a little bit further. Let $T$ denote the maximal torus of $G$ and $\{ P \}$ be its fixed point components. By using Atiyah-Bott localization formula, we can express the integral 

$$\int_M e^{\omega+i(\mu, \varphi)}$$
as the sum of the integrals over the fixed point components $\{P\}$. Here note that $\omega+i(\mu, \varphi)$ is a $G$-equivariant cohomology class. The interested reader can also interprete the derivation of the general formula in [Liu] by using the idea in this section.

\section{ Moduli spaces of flat connections}

In this section we will derive the symplectic volume formula and some results about the intersection numbers of the moduli space of flat connections on a principal flat $G$-bundle on a Riemann surface.  Some of the discussions about moduli spaces in this section are basically contained in [Liu1] and [Liu2] where we considered maps like 

$$f: \ \ G^{2g}\times G^{n}\rightarrow G $$ with 

$$ \ f(x_1, y_1, \cdots,x_g, y_g; z_1,\cdots z_n)=\prod_{j=1}^g[x_j, y_j]\prod_{i=1}^n z_ic_i z_i^{-1}$$ for certain fixed generic points $c_1, \cdots, c_n\in  G$. Here $G$ is a semisimple simply connected compact Lie group. We refer the reader to [Liu1] and [Liu2] for the details of the study of these maps by using heat kernels. Here let us consider an equivalent map given by 

$$f: \ \ G^{2g} \times \prod_{i=1}^n O_{c_i}\rightarrow G$$ with 
$$f(x_1, \cdots, y_g; z_1, \cdots, z_n)=\prod_{j=1}^g [x_j, y_j]\prod_{i=1}^n z_i$$
which is more commonly used to describe the moduli spaces. Here $x_j, \ y_j\in G$ and $z_i\in O_{c_i}$. Recall that $O_{c_i}$ is the conjugacy class of $c_i$ in $G$.

Now let us start to derive the formulas. We equip $G$ with the bi-invariant metric induced by the Killing form, then the explicit expression of the heat kernel on $G$ is given by  

$$H(t, x, y)=\frac{1}{|G|}\sum_{\lambda \in P_+}d_\lambda \cdot  \chi_\lambda(xy^{-1})e^{-tp_c(\lambda)}$$ where $|G|$ denotes the volume of $G$ and $P_+$ denotes all irreducible representations of $G$, which can be identified as a lattice in $t^*$, the dual of the Lie algebra of the maximal torus $T$. Also $d_\lambda$ and $\chi_\lambda$ denote the dimension and respectively the character of the representation $\lambda$, and  $p_c(\lambda)=|\lambda+\rho|^2-|\rho|^2$ with 

$$\rho=\frac{1}{2}\sum_{\alpha\in \Delta^+}\alpha $$ the half sum of the positive roots.

Recall that the moduli space we are interested in is just 

$${\cal M}_c={\cal M}_{c_1, \cdots, c_n}= f^{-1}(e)/G$$where $G$ acts on $G^{2g}\times \prod_{i=1}^n O_{c_i}$ by the conjugation $\gamma$:

$$\gamma: \ G\rightarrow G^{2g}\times \prod_iO_{c_i}$$ with 

$$\gamma(w)(x_1, \cdots, y_g; z_1, \cdots, z_n)=(wx_1w^{-1}, \cdots, wy_gw^{-1}; wz_1w^{-1},\cdots, wz_nw^{-1}).$$

Similarly we consider the integral 

$$I(t)=\int_{h\in G^{2g}\times \prod_jO_{c_j}}H(t, f(h), e) dh$$where $dh$ denotes the volume form of the bi-invariant metric on $G^{2g}\times \prod_jO_{c_j}$ induced from the Killing form. We will again perform the computation of $I(t)$ in two different ways: local and global.

(i) Local computation. Let $Z(G)$ denote the center of $G$ and $|Z(G)|$ denote the number of elements in $Z(G)$. As $t$ goes to $0$, a computation of the Gaussian integral as in [Liu1] gives us 

$$I(t)=\frac{|G|}{|Z(G)|}\int_{{\cal M}_{c}}d\nu_c+ O(e^{-\delta^2/4t})$$ where $d\nu_c$ is the Reidemeister torsion $\tau({\cal C}'_c)$ of the complex 

$${\cal C}'_c: \hspace{.3in} 0\rightarrow g\stackrel{d\gamma}{\rightarrow}g^{2g}\times \prod_j { Z}^\perp_{c_j}\stackrel{df}{\rightarrow}g \rightarrow 0$$where ${Z}^\perp_{c_j}$ denotes the tangent space to $O_{c_j}$ at $c_j$. It is clear that ${Z}^\perp_{c_j}$ can be viewed as the image of the map $$(I-\mbox{Ad}(c_j)):\ g\rightarrow g .$$ 

Here we take generic $c_j$'s, then $Z_{c_j}\simeq t$ and the tangent bundle $TM_c$ of $M_c$ is isomorphic to the first cohomology group of $\cal{C}_c'$. By using the Poincare duality of the cochain complexes of the Riemann surface, one can show that the torsion $\tau({\cal C}'_c)$ is related to the $L^2$-volume of ${\cal M}_c$ in the following way:

$$ \tau({\cal C}'_c)^2 =j(c)^2\|{\mbox{det}}\,T{\cal M}_c\|_{L^2}^2 $$ where 

$$j(c)^2=\prod_j|{ \mbox{ det}}(I-\mbox{Ad}( c_j))|$$can be considered as the torsion of the boundary of the Riemann surface. Here the determinant is taken by restricting to $Z^\perp_{c_j}$.

But from the definition we know that, up to a normalization by a factor of $2\pi$, the $L^2$-volume is exactly the symplectic volume:

 $$(2\pi)^{2N_c}\frac{{\omega_c^{N_c}}}{{N_c!}}=\|{\mbox{det}}\,T{\cal M}_c\|_{L^2} $$ 
for the moduli space ${\cal M}_c$. Here $N_c$ denotes the complex dimension of ${\cal M}_c$, and $\omega_c$ denotes the canonical symplectic form on ${\cal M}_c$ induced from the Poincare duality. So we have 

$$\tau({\cal C}'_c)=d\nu_c=(2\pi)^{2N_c}|j(c)|\frac{\omega^{N_c}}{N_c!}$$where 

$$|j(c)|= \prod_j|{ \mbox{ det}}(I-\mbox{Ad}( c_j))|^{\frac{1}{2}}.$$ We refer the reader to [W], [BL], [Liu1], [Liu2] for a proof of this relation. 

(ii) Global computation. The global computation is achieved by using the character relations:

$$\int_G\chi_\lambda(wyzy^{-1}z^{-1})dz=\frac{|G|}{d_\lambda}\chi_\lambda(wy)\chi_\lambda(y^{-1}),$$

$$\int_G\chi_\lambda(wy)\chi_\lambda(y^{-1})dy=\frac{|G|}{d_\lambda}\chi_\lambda(w)$$ and the formula

$$\int_{O_{c_j}}h(g)dv_g =\frac{|J(c_j)|}{|Z_{c_j}|}\int_G h(gc_j g^{-1})dg$$ for any continuous function $h$ on $O_{c_j}$. Here $Z_{c_j}$ denotes the stabilizer of $c_j$, $|Z_{c_j}|$ its induced Riemannian volume, and $dv_g$ is the induced volume form on $O_{c_j}$. And note that we have used the notation:

$$J(c_j)={ \mbox{ det}}(I-\mbox{Ad}(c_j)).$$

By applying these formulas inductively we get 

$$\int_{G^{2g}\times \prod_i O_{c_i}}\chi_\lambda(\prod_j [x_j, y_j]\prod_i z_i)\prod_j dx_j dy_j \prod_i dz_i$$

$$= |G|^{2g+n}\frac{ j(c)^2}{\prod_i|Z_{c_i}|}\frac{\prod_i \chi_\lambda(c_i)}{d_\lambda^{2g+n-1}}.$$

By putting the above computations together, we obtain the following formula which we also summarize as a proposition.

\vspace{.2in}

{\bf Proposition 2:} {\em We have the following formula for the symplectic volume of ${\cal M}_c$:}

$$\int_{{\cal M}_c}e^{\omega_c}=|Z(G)|\frac{|G|^{2g+n-2}|j(c)|}{(2\pi)^{2N_c}\prod_j |Z_{c_j}|}\sum_{\lambda\in P_+}\frac{\prod_j \chi_\lambda(c_j)}{d_\lambda^{2g+n-2}}e^{-tp_c(\lambda)}+O(e^{-\delta^2/4t}).$$

\vspace{.1in}

As in [Liu1] or [Liu2], we can then take derivatives with respect to the $c_j$'s on both sides of the above identity to get intersection numbers on the moduli space ${\cal M}_c$ or $ {\cal M}_u$ for $u$ an element in the center $Z(G)$. We refer the reader to [Liu1] and [Liu2] for the details of the derivation. Here we only mention an interesting vanishing theorem that we can easily obtain from the above formula. For simplicity let us take $n=1$. The general case is the same.

In fact let us introduce a function on the dual Lie algebra $t^*$ of the maximal torus $T$,  $\pi(\lambda)=\prod_{\alpha\in \Delta^+}<\alpha, \lambda>$, from which we construct a differential operator $ \pi(\partial)$ such that, for $C\in t$ with $u\,{\mbox{exp}}\,C=c $, 

$$\pi(\partial) e^{(\lambda,\ C)}=\pi(\lambda) e^{(\lambda,\ C)}$$ where $(\lambda, C)$ denotes the natural pairing. By applying $\pi(\partial)^{2g}$ to $I(t)$, we get, up to a constant,  

$$\pi(\partial)^{2g}I(t)=|Z(G)|\frac{|G|^{2g-1}|j(c)|}{(2\pi)^{2N_c} |Z_{c}|}\sum_{\lambda\in P_+}d_\lambda\chi_\lambda(c)e^{-tp_c(\lambda)}+O(e^{-\delta^2/4t}).$$

Note that, if $c\neq e$ where $e$ denotes the identity element of $G$, then the sum on the right hand side, which is just the heat kernel $H(t, c, e)$, has limit $0$ as $t$ goes to $0$. This in particular implies that ${\mbox{ lim}}_{t\rightarrow 0}I(t)$ is a piecewise polynomial in $C$ of degree at most $2g|\Delta^+|$. From this one easily deduces certain new vanishing theorems for the intersection numbers of the moduli spaces. See [Liu1] and [Liu2] for the details on the explicit formula for certain intersection numbers.

We remark that, to get more complete information about integrals on the moduli space ${\cal M}_c$, we can consider the integral 

$$I(t)=\int_{G^{2g}\times \prod_j O_{c_j}} F(h)H(t, f(h), e)dh$$ where $F(h)$ is some $G$-invariant function on $G^{2g}\times \prod_j O_{c_j}$. 

On one hand, as $t$ goes to $0$, $I(t)$ has limit given by  

$$\frac{|G|}{|Z(G)|} \int_{{\cal M}_c} \bar{F}e^{\omega_c}$$where $\bar{F}$ denotes the function on ${\cal M}_c$ induced by $F(h)$. This is the result of the local computation. Note that in principle the integral of any cohomology class on ${\cal M}_c$ can be written in the above form.

On the other hand, by Peter-Weyl theorem, we know that $F(h)$ can be expressed as the combinations of the characters of $G$. By using the product formula for characters

$$\chi_\mu\cdot \chi_\lambda=\sum_\nu C^\nu_{\mu \lambda}\cdot \chi_\nu$$where $C^\nu_{\mu \lambda}$ is the Clebsch-Gordon coefficients, we can perform the global computation to get an explicit infinite sum over $P_+$. As an example the interested reader may try to derive the formula by taking $F(h)=\prod_{j=1}^{k}\chi_{\mu_j}(x_j^{-1})$ with $k\leq g-1$.

\section{The push-forward of measures}

Now motivated by a conjecture of Diaconis, we consider the push-forward measure by the following map:

$$f: \ \ G^{2g}\rightarrow G, \ f(x_1, \cdots, y_g)=\prod_j[x_j,y_j].$$
Let $dh$ denote the Riemannian volume of $G^{2g}$, then Diaconis conjectures that the push-forward measure is given by the following formula:

$$f_*dh(x)= |G|^{2g-1}\sum_{\lambda\in P_+}\frac{\chi_\lambda(x^{-1})}{d_\lambda^{2g-1}}dx$$if the right hand side converges.  Here $dx$ denotes the biinvariant measure on $G$. Otherwise we should use the normalized limit

$$f_*dh(x)= |G|^{2g-1}{\mbox{lim}_{t\rightarrow 0}}\sum_{\lambda\in P_+}\frac{\chi_\lambda(x^{-1})}{d_\lambda^{2g-1}}e^{-tp_c(\lambda)}dx.$$

This conjecture can be proved by using the above heat kernel idea. Assume $f_*dh(x) =F(x)dx$ as forms on $G$,  then we only need to show that 

$$F(x) = |G|^{2g-1}\sum_{\lambda\in P_+}\frac{\chi_\lambda(x^{-1})}{d_\lambda^{2g-1}}.$$

It is clear that we have the identity

$$I(t)=\int_{G^{2g}}H(t, f(h), x)dh= \int_GH(t, y, x)F(y)dy.$$

Note that, as $t$ goes to $0$, the right hand side has limit given by $F(x)$, and the left hand side can be calculated by using the character relations. This immediately gives the above conjectured equality.

In fact we can prove similar formula for the push-forward of the Riemannian measure by the more general map 

$$f: \ \ G^{2g} \times \prod_j O_{c_j}\rightarrow G, \ \ f(x_1, \cdots, y_g; z_1, \cdots, z_n)=\prod_j [x_j, y_j]\prod_j z_j,$$ from which we get the following formula.

\vspace{.2in}

{\bf Proposition 3: } {\em Let $dh$ denotes the Riemannian measure on $  G^{2g} \times \prod_j O_{c_j}$. Then the following formula holds on $G$}

$$f_*dh(x)=|G|^{2g+n-1}\frac{j(c)^2}{\prod_j |Z_{c_j}|}\sum_{\lambda\in P_+}\frac{\prod_j \chi_\lambda(c_j)}{d_\lambda^{2g+n-1}}\chi_\lambda(x^{-1})dx . $$

\vspace{.1in}

Here we assume the right hand side converges, otherwise we should consider the normalized limit as above.

Now we want to consider a slightly different situation. Let $H_j, \ j=1, \cdots, n, $ be  subgroups of $G$, we consider the following map 

$$f: \ \ G^n \times \prod_j H_j\rightarrow G,\ \ f(x_1, \cdots, x_n; u_1, \cdots, u_n)=\prod_j x_j u_j x_j^{-1}.$$

The same argument of using heat kernel will give us an interesting formula for the push-forward measure. Indeed we consider the integral 

$$I(t)= \int_{h\in G^n \times \prod_j H_j}H(t, f(h), x)dh= \int_GH(t, y, x)F(y)dy.$$
 As $t$ goes to $0$, the right hand side is $F(x)$, and the left hand side can be calculated by using the character relations, from which we get 

\vspace{.2in}

{\bf Proposition 4:} {\em Let $dh$ denote the Riemaninn measure on $G^n \times \prod_j H_j$. Then the following formula holds:}

$$f_*dh(x)= {|G|^{n-1}}\sum_{\lambda\in P_+}\frac{\prod_j \int_{H_j}\chi_\lambda(u_j) du_j}{d_\lambda^{n-2}}\chi_\lambda(x^{-1})dx. $$ 

\vspace{.1in}

One may also consider the map 

$$f: \ G^{2g} \times G \rightarrow G$$ with  

$$f(x_1, \cdots, y_g; z)= \prod_j [x_j, y_j]z^2,$$ as well as the map 

$$f: \ G^{2g}\times G \times G\rightarrow G$$ with 

$$f(x_1, \cdots, y_g; w; z)= \prod_j [x_j y_j]wzw^{-1}z.$$

We leave these as exercises to the interested reader. In these cases one needs the formula 
$$\int_G \chi_\lambda(x^2)dx = |G|f_\lambda $$ where $f_\lambda =1$ if $\lambda$ has real structure, $-1$ if $\lambda$ has a quartenionic structure, $0$ otherwise. 

More generally, we can consider the problem of solving equations 
$$f_j (x_1, \cdots, x_m)=c_j, \ j=1, \cdots,n $$
in Lie groups. To understand the measure or number of solutions of these equations, we may apply the heat kernel on the corresponding group $G$ to the map

$$f: \ \ G^m \rightarrow G^n$$ with $f=(f_1, \cdots, f_n)$, and consider the integral 

$$I(t) = \int_{G^n}\prod_{j=1}^n H(t, f_j(h), c_j)dh.$$ In many cases we can derive various interesting formulas. As an example, let us consider the map 

$$f^{(n)} : \ \ G^n \rightarrow G, \ \ f^{(n)} (x_1, \cdots, x_n)=[x_1, [x_2, [\cdots, x_n]]]. $$ This is the $n$-commutator map. We will derive an inductive formula for the push-forward measure, Let us define $Q_n(x)$ by the identity :
$$f^{(n)}_*dh(x) =Q_n(x)dx.$$ Then we will prove  

\vspace{.2in} 

{\bf Proposition 5:} {\em The following formula holds:} 

$$Q_n(w)={\mbox{ lim}}_{t\rightarrow 0}\sum_{\lambda\in P_+}e^{-tp_c(\lambda)}\int_G\chi_\lambda(g)\chi_\lambda(g^{-1}w^{-1})Q_{n-1}(g) dg.$$

\vspace{.1in}
 
To prove this formula, we consider the integral 

$$I(t) = \int_{G^n}H(t, f^{(n)}(h), x)dh= \int_GH(t, y, x)Q_n(y)dy.$$ As $t$ goes to $0$, the right hand side is nothing but $Q_n(x)$, while the left hand side is reduced to the computation of the integrals like

$$ \int_{G^n} \chi_\lambda([x_1, [\cdots, x_n]]x^{-1})\prod_{j=1}^n dx_j$$ which, after integration with respect to $x_1$ by using the character relation,  is equal to 

$$\frac{|G|}{d_\lambda}\int_{G^{n-1}}\chi_\lambda([x_2, [\cdots, x_n]])\chi_\lambda([x_2, [\cdots, x_n]]^{-1}x^{-1})\prod_{j=2}^{n} dx_j$$

$$=\frac{|G|}{d_\lambda}\int_{G}\chi_\lambda(w)\chi_\lambda(w^{-1}x^{-1})f^{(n-1)}_*[\prod_{j=2}^{n} dx_j](w)$$

$$ =\frac{|G|}{d_\lambda}\int_{G}\chi_\lambda(w)\chi_\lambda(w^{-1}x^{-1})Q_{n-1}(w)dw.$$ Here the second identity is from the change of variable 

$$w= [x_2, [\cdots, x_n]].$$

By putting all of the above formulas together, we get the wanted formula.

\section{Counting solutions in finite groups}

It turns out that all of the results for compact Lie groups in the above section have analogues for  finite groups. Even the proofs are basically the same, if we replace the integrals in the last section by sums. 

As pointed to me by P. Diaconis, the heat kernel for a finite group $G$ is given by a similar expression:

$$H(t, x, y)= \frac{1}{|G|}\sum_{\lambda\in P_+}d_\lambda \cdot \chi_\lambda(xy^{-1}) e^{-tp_c(\lambda)}$$ where $|G|$ denotes the number of elements in $G$, $P_+$ still denotes all of the irreducible representations which is a finite set, and $p_c(\lambda)$ is a function on $P_+$. In fact the only property we need for this function is that its limit, as $t$ goes to $0$, is the delta function:

$$\delta(xy^{-1})=\frac{1}{|G|}\sum_{\lambda\in P_+}d_\lambda \cdot \chi_\lambda(xy^{-1})$$where $\delta(xy^{-1})=1$ if $x=y$ and $0$ otherwise. 

Let $O_{c_j}$ with $c_1, \cdots, c_n \in G$ denote the conjugacy class of $c_j$ in $G$. Let $S_{g, n}$ denote the number of solutions in $G^{2g}\times \prod_j O_{c_j}$of the equation 

$$\prod_{j=1}^g [x_j, y_j]\prod_{j=1}^n z_j =e$$ where $x_j, y_j \in G$ and $z_j \in O_{c_j}$.

This is related to the Chern-Simons theory with finite gauge group. Let $S$ be a compact Riemann surface. For $n=0$, $S_{g, n}$ counts the number of elements in the set ${\mbox{ Hom}}\, (\pi_1(S), G)$. We will derive a general formula for $S_{g, n}$. The special case when $n=0$ was proved in [FQ] by using topological field theory. As pointed out in [FQ], such formula was actually known to Serre.  

It is interesting to consider some variations of the above equation. Given subgroups $H_1, \cdots, H_n$ in $G$, we may consider the equation in $G^{n}\times \prod_j H_j$:

$$\prod_{j=1}^n x_jw_jx_j^{-1} =e$$ where $x_j \in G$ and $w_j \in H_j$. 

We can also consider the $n$-commutator equation 
$$[ x_1,[x_2, \cdots, x_n]]=e$$. 

In all of these cases we will give general formulas. Certainly one may try to find more general equations by combining the above equations together, or try to figure out equations of other types.

As in the compact Lie group cases, all of the formulas will follow naturally from the heat kernel trick. For simplicity and compatibility with last section, let us introduce a notation: for a finite set $G$ and any function $f$ on it, let us write 

$$\int_G f(g) dg = \sum_{g\in G} f(g).$$

We start from the first problem. We consider the map 

$$f: \ \ G^{2g}\times \prod_j O_{c_j} \rightarrow G$$
with 
$$ f(x_1, y_1, \cdots, x_g, y_g; z_1, \cdots, z_n)= \prod_{j=1}^n [x_j, y_j]\prod_{j=1}^n z_j,$$ 
and introduce the integral, more precisely the sum, 

$$I(t) = \int_{G^{2g}\times \prod_j O_{c_j} }H(t, f(h), e)dh.$$

(i) Local computation. As $t$ goes to $0$, the delta-function property of the heat kernel tells us that the limit is exactly the number $S_{g, n}$.

(ii) Global computation. On the other hand we can explicitly calculate the integral by using the formulas

$$\int_G\chi_\lambda(wyzyz^{-1})dz=\frac{|G|}{d_\lambda}\chi_\lambda(wy)\chi_\lambda(y^{-1}),$$

$$\int_G\chi_\lambda(wy)\chi_\lambda(y^{-1})dy=\frac{|G|}{d_\lambda}\chi_\lambda(w),$$and if $Z_{c_j}$ denotes the stabilizer of $c_j$, 

$$\int_{O_{c_j}}h(g)dv_g =\frac{1}{|Z_{c_j}|}\int_G h(gc_j g^{-1})dg.$$
Here as our convention, the integral means taking sums over $G$, and $|G|$, $|Z_{c_j}|$ denote the number of elements in $G$ and $Z_{c_j}$ respectively.

By putting the above two computations together, we get the following:

\vspace{.2in}

{\bf Proposition 6:} {\em We have the following formula:}

$$S_{g, n}= \frac{|G|^{2g+n-1}}{\prod_{j=1}^n |Z_{c_j}|}\sum_{\lambda\in P_+}\frac{\prod_{j=1}^n \chi_\lambda(c_j)}{d_\lambda^{2g+n-2}}.$$

\vspace{.1in}

For the second equation, let us define $Q_n(w)$ to be the number of solutions in $G^n$ for the $n$-commutator equation, then we will derive an induction formula which was first derived by Diaconis by using combinatorics: 

\vspace{.2in}

{\bf Proposition 7:} {\em The following formula holds: }

$$Q_n(w)=\sum_{\lambda\in P_+}\sum_{g\in G}\chi_\lambda(g)\chi_\lambda(g^{-1}w^{-1})Q_{n-1}(g).$$

\vspace{.1in}

To prove this we consider the map 

$$f^{(n)}: \ \ G^n \rightarrow G, \ \ f(x_1, \cdots, x_n)=[ x_1, [x_2,\cdots, x_n]]$$
and the function: 

$$I(t)=\int_{G^n} H(t, f^{(n)}(h), w)dh=\int_GH(t, y, w) Q_n(y)dy$$

(i) Local computation. As $t$ goes to $0$, the limit is clearly the function $Q_n(w)$.

(ii) Global computation. By using the character formulas, similar to the compact Lie group case, we get

$$\int_{G^n} \chi_\lambda(f(h)w^{-1})dh= \frac{|G|}{d_\lambda}\int_G\chi_\lambda(g)\chi_\lambda(g^{-1}w^{-1})Q_{n-1}(g)dg.$$
We then obtain the wanted formula by identifying the above two computations.

For the third equation, let $R_n$ be the number of solutions in $G^n \times \prod_{j=1}^nH_j$ such that, for $x_j \in G$ and $z_j \in H_j$,  
$$\prod_{j=1}^n x_j z_j x_j^{-1}=e.$$

Let us consider the map 

$$f: \ \ G^n \times \prod_{j=1}^nH_j \rightarrow G$$with 

$$ \ \ f(x_1, \cdots, x_n; z_1, \cdots, z_n)= \prod_{j=1}^n x_j z_j x_j^{-1}.$$

Again we consider the function

$$I(t) = \int_{ G^n \times \prod_{j=1}^nH_j} H(t, f(h), e)dh,$$ and perform the local and global computations.

(i) Local computation. As $t$ goes to $0$, the limit is given by the $R_n$.

(ii) Global computation. Still by using the character formulas we get the expression:

$$\int_{G^n\times\prod_i H_i}\chi_\lambda (\prod_{j=1}^n x_j z_j x_j^{-1})\prod_j dx_j dz_j= {|G|^n}\sum_{\lambda\in P_+}\frac{\prod_j\int_{H_j}\chi_\lambda( z_j)dz_j}{d_\lambda^{n-1}}.$$

So we have obtained our next proposition:

\vspace{.2in}

{\bf Proposition 8:} {\em We have the following expression for $R_n$:}

$$R_n= |G|^{n-1}\sum_{\lambda\in P_+}\frac{\prod_j\sum_{z_j\in H_j}\chi_\lambda(z_j)}{d_\lambda^{n-2}}.$$

\vspace{.1in}

More generally, we can find the number of solutions in finite groups of some general equations like 

$$f_j (x_1, \cdots, x_n) =c_j, \ \ c_, \ x_j \in G, \ j= 1, \cdots, m$$ by considering the map 

$$f: \ \ G^n \rightarrow G^m, \ f(h)=(f_1(h), \cdots, f_m(h))$$ and the function from the heat kernel:

$$I(t)= \int_{G^n} \prod_{j=1}^n H(t, f_j(h), c_j^{-1})dh.$$

Then similarly we can perform the global and local computations as above to get explicit formulas. 

We leave to the reader as an exercise to find the numbers of solutions in the equations like 

$$\prod_{j=1}^g [x_j, y_j] z^2= e, \ {\mbox{in }}\  G^{2g} \times G$$ and 

$$ \prod_{j=1}^g [x_j, y_j] wzw^{-1}z=e, \ {\mbox{in }}\ G^{2g}\times G\times G.$$

Another interesting problem is to consider the integrals, or sums, like  

$$I(t)=\int_{G^{2g}} F(h) H(t, f(h), e)dh$$for some $G$-invariant function $F$ on $G^{2g}$. As $t\rightarrow 0$, the limit of $I(t)$ compute 

$$\sum_{h\in {\cal M}} \bar{F}(h)$$ where ${\cal M}$ denotes the set of solutions of one of the equations in this section. By expressing $F(h)$ in terms of the characters of $G$, we can obtain certain explicit expression of this sum in terms of the representation data of $G$.

Note that the numbers of solutions of some of the above equations correspond to the Hurwitz numbers, the numbers of coverings of a given Riemann surface. In [St], the formula in Proposition 6 was used to derive the general case of the Brauer $p$-block conjecture. 

\section{Conluding remarks}

There are some other interesting cases to which one may consider to apply the above  heat kernel method. In [Liu2] we tried to find invariants for knots and $3$-manifolds. It may be interesting to combine this method with the work in [Hu]. 

As one can see that, symplectic structure has played the important role in our computations, such as in the derivation of the nonabelian localization formula and the intersection number formulas for moduli spaces. Recently P. Xu has explained to me some very general constructions of moment maps to Poisson manifolds. It should be interesting to apply our method to such general situation. An interesting example is to consider the group-valued moment map as introduced in [AMM], one may consider the integral of the pull-back of the  heat kernel on the Lie group $G$ by this map against the measure introduced in [AMM]. This should give us some quite interesting formulas relating the symplectic reduction and the representations of the Lie group. More precisely thhis should give us formulas similar to Witten's nonabelian localization formula. I thank Prof. S. Wu for the discussions concerning this point.

\end{document}